\documentclass{article}

\usepackage{amsthm}
\usepackage{amssymb}
\usepackage{amsxtra}

\newtheorem{thm}{Theorem}[section]
\newtheorem{lemma}{Lemma}[section]

\title{Higher moments of primes in short intervals I}
\author{Tsz Ho Chan}

\begin{document}
\maketitle
\begin{abstract}
In this article, we prove an ``equivalence'' between two higher
even moments of primes in short intervals under Riemann
Hypothesis. We also provide numerical evidence in support of these
asymptotic formulas.
\end{abstract}
\section{Introduction}

Recently, Montgomery and Soundararajan [\ref{MS2}] studied the
moments
$$M_k(N;h) := \sum_{n=1}^{N} (\psi(n+h) - \psi(n) - h)^k$$
where $k$ is a positive integer, $\psi(x) = \sum_{n \leq x}
\Lambda(n)$ and $\Lambda(n)$ is von Mangoldt lambda function. They
proved that, under a strong form of Hardy-Littlewood prime-$k$
tuple conjecture, for small $\epsilon > 0$, there is a $\delta >
0$,
\begin{equation}
\label{1.1} M_k(N;h) = \mu_k h^{k/2} \int_{1}^{N}
(\log{\frac{x}{h}} + B)^{k/2} dx + O_k(h^{k/2} N^{1-\epsilon})
\end{equation}
uniformly for $(\log N)^{15 k^2} \leq h \leq N^{1/k - \delta}$
where $\mu_k = 1 \cdot 3 \cdot \cdot \cdot (k-1)$ if $k$ is even,
and $\mu_k = 0$ if $k$ is odd. Here $B = 1 - C_0 - \log 2\pi$ and
$C_0$ denotes Euler's constant. One further expects that
(\ref{1.1}) holds uniformly for $N^\delta \leq h \leq
N^{1-\delta}$. This implies that, for $0 \leq x \leq N$, the
distribution of $\psi(x+h) -\psi(x)$ is approximately normal with
mean $h$ and variance $h \log{N/h}$. It contradicts with the
prediction of Cram\'{e}r's model of variance $h \log N$. In the
last section, we will show numerical evidence in support of
(\ref{1.1}).

Now, $M_k(X;h)$ can be written as
\begin{equation}
\label{1.2} \int_{1}^{X} (\psi(x+h) - \psi(x) - h)^k dx.
\end{equation}
We also consider the following moments:
\begin{equation}
\label{1.3} \widetilde{M}_k(X;\delta) := \int_{1}^{X} (\psi(x +
\delta x) - \psi(x) - \delta x)^k dx.
\end{equation}
Goldston and Montgomery [\ref{GM}] showed that, under Riemann
Hypothesis (RH), the stronger form of the Pair Correlation
Conjecture as formulated by Montgomery [\ref{M}] is equivalent to
an asymptotic formula for (\ref{1.2}) in $X^\epsilon \leq h \leq
X^{1-\epsilon}$ or an asymptotic formula for (\ref{1.3}) in
$X^{-1+\epsilon} \leq \delta \leq X^{-\epsilon}$ when $k=2$. The
author generalized these to include the second main terms in
[\ref{C}] (again only when $k=2$). So, the main purpose of this
paper is to prove the ``equivalence'' between an asymptotic
formula for (\ref{1.2}) and an asymptotic formula for (\ref{1.3})
in appropriate ranges of $h$ and $\delta$ for any positive even
integer $k$. Roughly speaking, we have
\begin{thm}
\label{theorem1.1} Let $k$ be a positive even integer. Assuming
RH, the following are equivalent:
$$(i) \; \int_{1}^{X} (\psi(x+h) - \psi(x) - h)^k dx \sim \mu_k
h^{k/2+1} \int_{E}^{X/h} (\log{\frac{x}{E}})^{k/2} dx$$ holds
uniformly for $X^\epsilon \leq h \leq X^{1-\epsilon}$.
$$(ii) \; \int_{1}^{X} (\psi(x+\delta x) - \psi(x) - \delta x)^k
dx \sim \frac{\mu_k}{\frac{k}{2}+1} X^{k/2+1} \delta^{k/2}
\Bigl(\log{\frac{1}{E \delta}}\Bigr)^{k/2}$$ holds uniformly for
$X^{-1+\epsilon} \leq \delta \leq X^{-\epsilon}$.
\end{thm}
Here $E = 2\pi e^{C_0 - 1}$. Our method of proof replaces the
brute-force calculations in [\ref{C}]. We will assume RH
throughout this paper and $k$ being a positive even integer unless
stated otherwise.

This work is part of the author's $2002$ PhD thesis with some
improvements.
\section{Some preparations}

First of all, $\psi(x) = x + O(x^{1/2} \log^2{x})$ by RH (see
[\ref{vK}]). One has the following:
\begin{equation}
\label{ineq1} \int_{1}^{X} (\psi(x+\delta x) - \psi(x) - \delta x
)^k dx \ll X^{{k/2}+1} \log^{2k}{X}
\end{equation}
for $0 \leq \delta \leq 1$, and
\begin{equation}
\label{ineq2} \int_{1}^{X} (\psi(x+h) - \psi(x) - h)^k dx \ll
X^{{k/2}+1} \log^{2k}{X}
\end{equation}
for $0 \leq h \leq X$. Also, estimating trivially, we have
\begin{equation}
\label{ineq3} \int_{1}^{X} (\psi(x + \delta x) - \psi(x) - \delta
x)^k dx \ll \int_{1}^{X} (\delta x \log{X})^k dx \ll \delta^k
X^{k+1} \log^k{X}
\end{equation}
for $0 \leq \delta \leq 1$, and
\begin{equation}
\label{ineq4} \int_{1}^{X} (\psi(x + h) - \psi(x) - h)^k dx \ll
\int_{1}^{X} (h \log{X})^k dx \ll h^k X \log^k{X}
\end{equation}
for $0 \leq h \leq X$. We also need some lemmas.
\begin{lemma}
\label{lemma2.1} For any differentiable function $f(u)$ and any $0
\leq \eta \leq 1$,
$$\int_{T}^{(1+\eta)T} f(u) du = \eta T f(T) + O\bigl(\eta^2 T^2 \mathop{max}_{T
\leq t \leq (1+\eta)T} |f'(t)|\bigr).$$
\end{lemma}

Proof: By mean-value theorem, we have
\begin{eqnarray*}
\int_{T}^{(1+\eta)T} f(u) du &=& \eta T f(T + \xi T) \mbox{ where
} 0 \leq \xi \leq \eta \\
&=& \eta T \bigl(f(T) + \xi T f'(T + \xi' T)\bigr) \mbox{ where }
0 \leq \xi' \leq \xi,
\end{eqnarray*}
and the lemma follows.
\begin{lemma}
\label{lemma2.2} For any positive integer $k$, we have
$$x^k-y^k=(x-y)P(x,y) + (x-y)^k$$
where $P(x,y)$ is some homogeneous polynomial of degree $k-1$.
\end{lemma}

Proof: By Factor Theorem, $z-1$ divides $z^k -1 -(z-1)^k$. So,
$$(z-1)P(z) = z^k -1 -(z-1)^k$$
for some integer polynomial $P(z)$ of degree $k-1$. Set $z = {x
\over y}$ and multiply both sides by $y^k$, we get the desired
result.
\begin{lemma}
\label{lemma2.3} For any positive integer $k$, and any
non-negative real numbers $\alpha$, $a_1,a_2,...,a_k$, we have
$$(a_1+a_2+...+a_k)^{\alpha} {\ll}_{\alpha ,k} \; a_1^{\alpha} + a_2^{\alpha}
+...+ a_k^{\alpha}.$$ Here, ${\ll}_{\alpha ,k}$ means that the
implicit constant may depend on $\alpha$ and $k$ but not on any
$a_i$'s.
\end{lemma}

Proof: Without loss of generality, suppose that $a_1$ is the
largest among the $a_i$'s. Then
$$(a_1+a_2+...+a_k)^{\alpha} \leq (k a_1)^{\alpha} \leq k^{\alpha}
(a_1^{\alpha}+a_2^{\alpha}+ ... +a_k^{\alpha}).$$
\section{(i) $\Rightarrow$ (ii)}

Throughout this and the next section, we think of $k$ as fixed.
\begin{thm}
\label{theorem3.1} Assume RH. If, for some small $\epsilon >
\epsilon_1 > 0$ (small in terms of $k$),
\begin{equation}
\label{3.1} \int_{1}^{X} (\psi(x+h)-\psi(x)-h)^k dx = \mu_k h^{k/2
+ 1} \int_{E}^{X/h} \bigr(\log{\frac{x}{E}}\bigr)^{k/2} dx +
O_k(h^{k/2} X^{1-\epsilon_1})
\end{equation}
holds uniformly for $X^{\epsilon} \leq h \leq X^{1-\epsilon}$,
then
\begin{equation}
\label{3.1result} \int_{1}^{X} (\psi(x + \delta x) - \psi(x) -
\delta x)^k dx = \frac{\mu_k}{\frac{k}{2}+1} X^{{k/2}+1}
\delta^{k/2} \Bigl( \log{\frac{1}{E\delta}} \Bigr) ^{k/2} +
O_k(\delta^{k/2} X^{{k/2}+1-\epsilon_2})
\end{equation}
holds uniformly for $X^{-1 +2\epsilon + 2\epsilon_1} \leq \delta
\leq X^{-\epsilon}/2$ with some $\epsilon_2 > 0$.
\end{thm}

Proof: Our method is that of Saffari and Vaughan [\ref{SV}]
employed in [\ref{GM}] and [\ref{C}]. Let $f(x,h)
=\psi(x+h)-\psi(x)-h$. Let $X^{-1+2\epsilon+\epsilon_1} \leq
\Delta \leq X^{-\epsilon}$. Say $\Delta = X^{-\mu}$ for some
$\epsilon \leq \mu \leq 1 - 2\epsilon - \epsilon_1$. We want to
calculate
\begin{equation}
\label{3.2} \int_{V/2}^{V} \int_{0}^{\Delta} (\psi(x+\delta x) -
\psi(x) - \delta x)^k d\delta \; dx.
\end{equation}
Substituting $h=\delta x$, (\ref{3.2}) becomes
$$\int_{\Delta V/2}^{\Delta V} \int_{h/\Delta}^{V} {f(x,h)^k \over x}
dx \; dh + \int_{0}^{\Delta V/2} \int_{V/2}^{V} {f(x,h)^k \over x}
dx \; dh$$
$$=\int_{\Delta V/2}^{\Delta V} \int_{h/\Delta}^{V} +
\int_{V^{\epsilon}}^{\Delta V/2} \int_{V/2}^{V} +
\int_{0}^{V^{\epsilon}} \int_{V/2}^{V} = I_1 + I_2 + I_3.$$ By
integration by parts, we have from (\ref{3.1}) that
\begin{eqnarray*}
& & \int_{U}^{V} {f(x,h)^k \over x} dx \\
&=& \Bigl[{1 \over x} \int_{1}^{x} f(u,h)^k du {\Bigr]}_{U}^{V} +
\int_{U}^{V}
\Bigl(\int_{1}^{x} f(u,h)^k du \Bigr){1 \over x^2} dx \\
&=& \mu_k h^{k/2+1} \Bigl[{1 \over V} \int_{E}^{V/h} \bigl(\log{x
\over E}\bigr)^
{k/2} dx-{1 \over U}\int_{E}^{U/h} \bigl(\log{x \over E}\bigr)^{k/2} dx\Bigr] \\
& &+ \mu_k h^{k/2+1}\int_{U}^{V} {1 \over x^2}
\int_{E}^{x/h}\bigl(\log{u \over E}
\bigr)^{k/2} du \; dx + O_{k}(U^{-\epsilon_1} h^{k/2}) \\
&=& T_1 + T_2 + O_{k}(U^{-\epsilon_1} h^{k/2})
\end{eqnarray*}
as long as $V^{\epsilon} \leq h \leq U^{1-\epsilon}$ with $U \leq
V \leq 2U$.
\begin{eqnarray*}
T_2 &=& \mu_k h^{k/2} \int_{U/h}^{V/h} \int_{E}^{y} \bigl(\log{u
\over E}\bigr)^
{k/2} du d({-1 \over y}) \\
&=& \mu_k h^{k/2} \Bigl[-{h \over V} \int_{E}^{V/h} \bigl(\log{u
\over E}\bigr)^
{k/2} du + {h \over U} \int_{E}^{U/h} \bigl(\log{u \over E}\bigr)^{k/2} du \\
& &+ \int_{U/h}^{V/h} {(\log{y/E})^{k/2} \over y} dy \Bigr].
\end{eqnarray*}
Therefore,
\begin{equation*}
\int_{U}^{V} {f(x,h)^k \over x} dx = {\mu_k \over {k \over 2}+1}
h^{k/2} \Bigl[\bigl(\log{V \over Eh}\bigr)^{k/2+1} -\bigl(\log{U
\over Eh}\bigr)^{k/2+1} \Bigr] + O_k(U^{-\epsilon_1} h^{k/2})
\end{equation*}
as long as $V^{\epsilon} \leq h \leq U^{1-\epsilon}$ with $U \leq
V \leq 2U$. Thus, for $I_1$ and $I_2$ to work, we need
\begin{equation}
\label{3cond} V^{\epsilon} \leq {\Delta V \over 2} \leq \Bigl({V
\over 2} \Bigr)^{1-\epsilon}, \mbox{ and } V^{\epsilon} \leq h
\leq \Bigl({h \over \Delta }\Bigr)^{1-\epsilon} \mbox{ for }
{\Delta V \over 2} \leq h \leq \Delta V.
\end{equation}
Since $\Delta = X^\mu$, for $X^{\mu + \epsilon} \leq V \leq X$,
one can check that (\ref{3cond}) are satisfied. Therefore,
\begin{eqnarray*}
I_1 &=& {\mu_k \over {k \over 2}+1} \Bigl[\int_{\Delta
V/2}^{\Delta V} h^{k/2} \bigl(\log{V \over Eh}\bigr)^{k/2+1} dh
-\int_{\Delta V/2}^{\Delta V} h^{k/2}
\bigl(\log{1 \over E\Delta} \bigr)^{k/2+1} dh \Bigr] \\
& &+ O_{k}({\Delta}^{k/2+1} V^{k/2+1-\epsilon_1}), \\
I_2 &=& {\mu_k \over {k \over 2}+1} \Bigl[\int_{0}^{\Delta V/2}
h^{k/2} \bigl(\log{V \over Eh}\bigr)^{k/2+1} dh -\int_{0}^{\Delta
V/2} h^{k/2} \bigl(
\log{V/2 \over Eh}\bigr)^{k/2+1} dh\Bigr] \\
& &+ O_{k}(V^{({k/2}+1)\epsilon} (\log{V})^{{k/2}+1}) +
O_{k}({\Delta}^{{k/2}
+1} V^{{k/2}+1-\epsilon_1}), \\
I_3 &\ll& V^{\epsilon} V^{k \epsilon} \log^k{V},
\end{eqnarray*}
Combining these, (\ref{3.2}) equals
$${\mu_k \over {k \over 2}+1} \biggl[ \int_{0}^{\Delta V} h^{k/2}\bigl(\log{V
\over Eh}\bigr)^{k/2+1} dh -\int_{0}^{\Delta V/2} h^{k/2}
\bigl(\log{V/2 \over Eh}\bigr)^{k/2+1} dh$$
$$- \bigl(\log{1 \over E\Delta}\bigr)^{k/2+1} \int_{\Delta V/2}^{\Delta V}
h^{k/2} dh \biggr] + O\bigl(V^{(k+1)\epsilon} \log^k {V} \bigr) +
O_{k} (\Delta^{k/2+1} V^{k/2+1-\epsilon_1})$$ for $X^{\mu +
\epsilon} \leq V \leq X$. Now, replacing $V$ by $X2^{-l}$ in the
above, summing over $0 \leq l \leq M=[{ (1 - \mu - \epsilon)
\log{X} \over \log{2}}]$,
\begin{equation}
\label{3.4}
\begin{split}
&\int_{X/2^M}^{X} \int_{0}^{\Delta} \bigl(\psi(x+\delta x)
-\psi(x) -\delta x \bigr)^k d\delta \; dx \\
=& {\mu_k \over {k \over 2}+1} \biggl[\int_{0}^{\Delta X} h^{k/2}
\bigl(\log{X \over Eh}\bigr)^{k/2+1} dh - {1 \over {k \over 2}+1}
\bigl(\log{1 \over E\Delta}
\bigr)^{k/2+1} (\Delta X)^{k/2+1} \biggr]\\
&+ O_k(X^{(k+1)\epsilon} \log^k {X}) + O_k({\Delta}^{k/2+1}
X^{k/2 + 1 - \epsilon_1}) \\
=& {\mu_k \over {k \over 2}+1} \int_{0}^{\Delta X} h^{k/2}
\bigl(\log{X \over Eh} \bigr)^{k/2} dh + O_k({\Delta}^{k/2+1}
X^{k/2+1-\epsilon_1})
\end{split}
\end{equation}
by integration by parts, and as $\Delta = X^{-\mu} \geq X^{-1 +
2\epsilon + \epsilon_1}$,
$$X^{(k+1)\epsilon} \log^k {X} \ll X^{(k+2)\epsilon + (k/2)
\epsilon_1} \ll \Delta^{k/2+1} X^{k/2+1-\epsilon_1}.$$ Using
(\ref{ineq3}),
\begin{equation}
\label{3.5} \int_{1}^{X/2^M} \int_{0}^{\Delta} (\psi(x+\delta x) -
\psi(x)- \delta x)^k d\delta \; dx \ll_k \Delta^{k+1} (X^{\mu +
\epsilon})^{k+1} \log^k {X}.
\end{equation}
But, since $\mu \leq 1 - 2\epsilon - \epsilon_1$,
\begin{equation}
\label{3.6}
\begin{split}
\Delta^{k/2+1} X^{k/2+1-\epsilon_1} &= \Delta^{k+1} X^{k/2 + 1 +
(k/2) \mu - \epsilon_1} = \Delta^{k+1} X^{(k+1)\mu + (k/2+1)
(1-\mu) - \epsilon_1} \\
&\geq \Delta^{k+1} X^{(k+1)\mu + (k/2+1) (2\epsilon + \epsilon_1)}
\gg \Delta^{k+1} (X^{\mu + \epsilon})^{k+1} \log^k {X}.
\end{split}
\end{equation}
Combining (\ref{3.4}), (\ref{3.5}) and (\ref{3.6}), we have
\begin{equation}
\label{3.7}
\begin{split}
& \int_{0}^{\Delta} \int_{1}^{X} (\psi(x+\delta x) - \psi(x)
-\delta x)^k dx \; d\delta \\
=& {\mu_k \over {k \over 2}+1} \int_{0}^{\Delta X} h^{k/2}
\bigl(\log{X \over Eh} \bigr)^{k/2} dh + O_k({\Delta}^{k/2+1}
X^{k/2+1-\epsilon_1})
\end{split}
\end{equation}
for $X^{-1 + 2\epsilon + \epsilon_1} \leq \Delta \leq
X^{-\epsilon}$.

\bigskip

We now deduce (\ref{3.1result}) from (\ref{3.7}). Set $\eta =
X^{-2 \epsilon_1 / 3}$. By Lemma \ref{lemma2.1}, one has for
$X^{-1 + 2\epsilon + 2\epsilon_1} \leq \Delta \leq
X^{-\epsilon}/2$,
\begin{equation}
\label{3.8}
\begin{split}
& \int_{\Delta}^{(1+\eta)\Delta} \int_{1}^{X} (\psi(x+\delta
x) - \psi(x) -\delta x)^k dx \; d\delta \\
=& {\mu_k \over {k \over 2}+1} \int_{\Delta X}^{(1+\eta)\Delta X}
h^{k/2} \bigl( log{X \over Eh} \bigr)^{k/2} dh +
O_k(\Delta^{k/2+1} X^{k/2+1-\epsilon_1}) \\
=& {\mu_k \over {k \over 2}+1} (\Delta X)^{{k/2}+1} \Bigl(\log{1
\over E\Delta} \Bigr)^{k/2} \eta + O_{k}\Bigl(\eta^2 (\Delta
X)^{k/2+1}(\log{1 \over \Delta})^{k/2} \Bigr) \\
&+ O_{k}(\Delta^{k/2+1} X^{k/2+1-\epsilon_1}).
\end{split}
\end{equation}
Let $g(x,\delta x)=f(x,\Delta x)$ for $\Delta \leq \delta \leq
(1+\eta)\Delta$. Then one can easily check that $f(x,\delta
x)-g(x,\delta x)=f\bigl((1+\Delta)x,(\delta-\Delta)x \bigr)$. So,
\begin{equation}
\label{3.9}
\begin{split}
& \int_{\Delta}^{(1+\eta)\Delta} \int_{1}^{X} \bigl(f(x,\delta
x)-g(x, \delta x ) \bigr)^k dx d\delta = \int_{0}^{\eta \Delta
\over 1+ \Delta} \int_{1+\Delta}^
{(1+\Delta)X} f(x,\delta x)^k dx d\delta \\
&\ll_{k} (\eta X \Delta)^{{k/2}+1} (\log{1 \over \eta
\Delta})^{k/2}
\end{split}
\end{equation}
by (\ref{3.7}), the choice of $\eta$ and the range of $\Delta$.
Thus, by Lemma \ref{lemma2.3}, (\ref{3.8}) and (\ref{3.9}),
\begin{equation}
\label{3.10}
\begin{split}
\int_{\Delta}^{(1+\eta)\Delta} \int_{1}^{X} g(x,\delta x)^k dx \;
d\delta &\ll_k
\int \int |f(x,\delta x)|^k + \int \int |f(x,\delta x) - g(x,\delta x)|^k \\
\ll_{k}& \eta X^{{k/2}+1} \Delta^{{k/2}+1} (\log{1 \over
\Delta})^{k/2} + \Delta^{{k/2}+1} X^{{k/2}+1 -\epsilon_1}.
\end{split}
\end{equation}
By Lemma \ref{lemma2.2} and Holder's inequality,
\begin{eqnarray*}
& &\int_{\Delta}^{(1+\eta)\Delta} \int_{1}^{X} f(x,\delta x)^k -
g(x,\delta x)^k dx \; d\delta \\
&=&\int_{\Delta}^{(1+\eta)\Delta} \int_{1}^{X} P(f,g) (f-g) +
\int_{\Delta}^{(1+\eta)\Delta} \int_{1}^{X} (f-g)^k \\
&\ll&\Bigl(\int \int |P(f,g)|^{k/(k-1)}\Bigr)^{(k-1)/k} \Bigl(\int
\int |f-g|^k
\Bigr)^{1/k} + \int \int |f-g|^k \\
&=& J_1^{(k-1)/k} J_2^{1/k} + J_2
\end{eqnarray*}
where $P(x,y)$ is a homogeneous polynomial of degree $k-1$.
$$J_2 \ll_{k} \eta^{{k/2}+1} X^{{k/2}+1} \Delta^{{k/2}+1} (\log{1 \over \eta
\Delta})^{k/2}$$ by (\ref{3.9}). And
\begin{eqnarray*}
J_1 &\ll_k& \int \int \bigl(\sum_{i+j=k-1} |f|^i |g|^j \bigr)^{k/(k-1)} \\
&\ll_k& \sum_{i+j=k-1} \int \int \bigl(|f| + |g|\bigr)^k \mbox{ by
binomial
theorem} \\
&\ll_k& \int \int f^k + g^k \mbox{ by Lemma \ref{lemma2.3}}\\
&\ll_k& \eta X^{{k/2}+1} \Delta^{{k/2}+1}(\log{1 \over \eta
\Delta})^{k/2} + \Delta^{{k/2}+1} X^{{k/2}+1 -\epsilon_1} \mbox{
by (\ref{3.8}) and (\ref{3.10})} \\
&\ll_k& \eta X^{{k/2}+1} \Delta^{{k/2}+1}(\log{1 \over \eta
\Delta})^{k/2} \mbox{ as } \eta=X^{-2\epsilon_1 /3}.
\end{eqnarray*}
Consequently, by Lemma \ref{lemma2.3},
\begin{equation}
\label{3.11} \int \int f^k-g^k \ll_k X^{{k/2}+1} \Delta^{{k/2}+1}
\eta^{3/2}(\log{1 \over \eta \Delta})^{k/2}.
\end{equation}
Therefore, by (\ref{3.11}) and (\ref{3.8}),
\begin{eqnarray*}
& &\eta \Delta \int_{1}^{X} (\psi(x+\Delta x) - \psi(x) -
\Delta x )^k dx = \int \int g^k \\
&=&\int \int f^k + O\Bigl(X^{{k/2}+1} \Delta^{{k/2}+1}
\eta^{3/2}(\log{1 \over \eta \Delta})^{k/2}\Bigr) \\
&=&{\mu_k \over {k\over 2}+1} (X \Delta)^{{k/2}+1} \Bigl(\log{1
\over E\Delta} \Bigr)^{k/2} \eta
+ O_k(\Delta^{k/2+1} X^{k/2+1-\epsilon_1})\\
& &+ O_k\Bigl(X^{k/2+1} \Delta^{k/2+1} \eta^{3/2} (\log{1 \over
\eta \Delta})^ {k/2}\Bigr).
\end{eqnarray*}
Dividing through by $\eta \Delta$, we have
\begin{equation*}
\begin{split}
&  \int_{1}^{X} (\psi(x+\Delta x) - \psi(x) - \Delta x)^k dx \\
=& {\mu_k \over {k\over 2}+1} X^{{k/2}+1} \Delta^{k/2}
\Bigl(\log{1 \over E \Delta}\Bigr)^{k/2} + O_k\Bigl({\Delta^{k/2}
X^{k/2 + 1 - \epsilon_1} \over \eta} \Bigr) \\
&+ O_k\Bigl(X^{k/2+1} \Delta^{k/2} \eta^{1/2} (\log{1 \over \eta
\Delta})^ {k/2}\Bigr).
\end{split}
\end{equation*}
Finally, recall $\eta = X^{-2\epsilon_1/3}$, one has the error
terms $\ll X^{{k/2}+1-{\epsilon_1/4}} \Delta^{k/2}$. So, the
theorem is true with $\epsilon_2 = {\epsilon_1 \over 4}$.
\section{(ii) $\Rightarrow$ (i)}

\begin{thm}
\label{theorem4.1} Assume RH. If, for some small $\epsilon
> \epsilon_1 > 0$ (small in terms of $k$),
\begin{equation}
\label{4.1} \int_{1}^{X} (\psi(x+\delta x) - \psi(x) - \delta x
)^k dx = {\mu_k \over {k\over 2}+1} X^{{k/2}+1} \delta^{k/2}
\Bigl( \log{1 \over E\delta}\Bigr)^{k/2} + O_k(\delta^{k/2}
X^{k/2+1-\epsilon_1})
\end{equation}
holds uniformly for $X^{-1+\epsilon} \leq \delta \leq
X^{-\epsilon}$, then
\begin{equation}
\label{4.2} \int_{1}^{X} (\psi(x+h)-\psi(x)-h)^k dx = \mu_k
h^{k/2+1} \int_{E}^ {X/h} \bigl(\log{x \over E}\bigr)^{k/2} dx +
O_k(h^{k/2} X^{1-\epsilon_2})
\end{equation}
holds uniformly for $X^{2\epsilon + \epsilon_1} \leq h \leq
X^{1-(k/2+1) \epsilon - 2\epsilon_1} / 2$ with some $\epsilon_2 >
0$.
\end{thm}

Proof: Let $f(x,h)=\psi(x+h)-\psi(x)-h$. Let $X^{2\epsilon} \leq H
\leq X^{1 - (k/2+1)\epsilon - 2\epsilon_1)}$. Say $H = X^\mu$ for
some $2\epsilon \leq \mu \leq 1 - (k/2+1)\epsilon - 2\epsilon_1$.
First, we calculate
\begin{equation}
\label{4.3} \int_{V/2}^{V} \int_{0}^{H} (\psi(x+h) - \psi(x) - h
)^k dh \; dx.
\end{equation}
Substituting $\delta = {h \over x}$, (\ref{4.3}) becomes
$$\int_{H/V}^{2H/V} \int_{V/2}^{H/\delta} f(x,\delta x)^k x dx \; d\delta +
\int_{0}^{H/V} \int_{V/2}^{V} f(x,\delta x)^k x dx \; d\delta$$
$$=\int_{H/V}^{2H/V} \int_{V/2}^{H/\delta} + \int_{({V \over 2})^{-1+\epsilon}}^
{H/V} \int_{V/2}^{V} + \int_{0}^{({V \over 2})^{-1+\epsilon}}
\int_{V/2}^{V} = I_1+I_2+I_3.$$ By integration by parts, we have
from (\ref{4.1}) that
\begin{eqnarray*}
& & \int_{U}^{V} f(x,\delta x)^k x dx \\
&=& {\mu_k \over {{k \over 2}+2}} \delta^{k/2} \Bigl(\log{1 \over
E\delta}\Bigr)^ {k/2} \Bigl[V^{{k/2}+2} - U^{{k/2}+2}\Bigr] +
O_k(\delta^{k/2} V^{{k/2} +2-\epsilon_1})
\end{eqnarray*}
as long as $U^{-1+\epsilon} \leq \delta \leq V^{-\epsilon}$ with
$U \leq V \leq 2U$. In order for this to work for $I_1$ and $I_2$,
we need
\begin{equation}
\label{4cond} \Bigl({V \over 2}\Bigr)^{-1+\epsilon} \leq {H \over
V} \leq V^{-\epsilon}, \mbox{ and } \Bigl({V \over
2}\Bigr)^{-1+\epsilon} \leq \delta \leq \Bigl({H \over
\delta}\Bigr)^{-\epsilon} \mbox{ for } {H \over V} \leq \delta
\leq {2H \over V}.
\end{equation}
Since $H = X^\mu$, one can check that (\ref{4cond}) are satisfied
for $X^{\mu + \epsilon} \leq V \leq X$. Thus,
\begin{eqnarray*}
I_1 &=& {\mu_k \over {k \over 2}+2} \int_{H/V}^{2H/V}
\Bigl[\bigl({H \over \delta}\bigr)^{{k/2}+2} - \bigl({V \over
2}\bigr)^{{k/2}+2} \Bigr]
\delta^{k/2} \Bigl(\log{1 \over E\delta}\Bigr)^{k/2} d\delta \\
& &+ O_k(H^{{k/2}+1} V^{1-\epsilon_1}), \\
I_2 &=& {\mu_k \over {k \over 2}+2} \Bigl[V^{{k/2}+2} - \bigl({V
\over 2}\bigr)^ {{k/2}+2} \Bigr] \int_{0}^{H/V} \delta^{k/2}
\Bigl(\log{1 \over E\delta}\Bigr)^{k/2} d\delta \\
& &+ O_k(V^{1+({k/2}+1) \epsilon} (\log{V})^{k/2}) +
O_k(H^{{k/2}+1} V^{1- \epsilon_1}), \\
I_3 &\ll& V^{-1+\epsilon} V^2 V^{k\epsilon} \log^k {V} =
V^{1+(k+1)\epsilon} \log^k {V}.
\end{eqnarray*}
Let $\nu_k = \mu_k / ({k \over 2}+2)$, then (\ref{4.3})
\begin{eqnarray*}
&=& \nu_k H^{{k/2}+2} \int_{H/V}^{2H/V} {1 \over \delta^2}
\Bigl(\log{1 \over E\delta}\Bigr)^{k/2} d\delta + \nu_k
V^{{k/2}+2} \int_{0}^{H/V}
\delta^{k/2} \Bigl(\log{1 \over E\delta}\Bigr)^{k/2} d\delta \\
& &- \nu_k \Bigl({V \over 2}\Bigr)^{{k/2}+2} \int_{0}^{2H/V}
\delta^{k/2} \Bigl(\log{1 \over E\delta}\Bigr)^{k/2} d\delta \\
& &+O_k(V^{1+(k+1)\epsilon} \log^k {V}) + O_k(H^{{k/2}+1}
V^{1-\epsilon_1})
\end{eqnarray*}
when $X^{\mu + \epsilon} \leq V \leq X$. Now, replacing $V$ by
$X2^{-l}$ in the above, summing over $0 \leq l \leq M =[{(1 - \mu
- \epsilon) \log{X} \over \log{2}}]$,
\begin{equation}
\label{4.4}
\begin{split}
& \int_{X/2^M}^{X} \int_{0}^{H} (\psi(x+h) - \psi(x) - h)^k dh \; dx \\
=& \nu_k H^{{k/2}+2} \int_{H/X}^{2^{M+1}H/X} {1 \over \delta^2}
\Bigl(\log{1 \over E\delta}\Bigr)^{k/2} d\delta + \nu_k
X^{{k/2}+2} \int_{0}^{H/X} \delta^{k/2} \Bigl(\log{1 \over E\delta}
\Bigr)^{k/2} d\delta \\
&- \nu_k \Bigl({X \over 2^{M+1}}\Bigr)^{{k/2}+2}
\int_{0}^{2^{M+1}H/X} \delta^{k/2} \Bigl(\log{1 \over
E\delta}\Bigr)^{k/2} d\delta + O_k(H^{{k/2}+1} X^{1-\epsilon_1}) \\
=& \nu_k H^{{k/2}+2} \int_{H/X}^{1/E} {1 \over \delta^2}
\Bigl(\log{1 \over E\delta}\Bigr)^{k/2} d\delta + \nu_k
X^{{k/2}+2} \int_{0}^{H/X} \delta^{k/2} \Bigl(\log{1 \over E\delta}
\Bigr)^{k/2} d\delta \\
&+ O_k(H^{{k/2}+1} X^{1-\epsilon_1}).
\end{split}
\end{equation}
because, as $X^{2\epsilon} \leq H$, $X^{1+(k+1)\epsilon} \log^k
{X} \ll H^{k/2+1} X^{1-\epsilon_1}$. Also, the terms involving
$2^{M+1}$ are absorbed into the error term as $\mu \leq 1 -
\epsilon - 2\epsilon_1$. Using (\ref{ineq2}),
\begin{equation}
\label{4.5} \int_{1}^{X/2^{M}} \int_{0}^{H} \bigl(\psi(x+h)
-\psi(x) -h \bigr)^k dh \; dx \ll H (X^{\mu + \epsilon})^{k/2+1}
\log^{2k} {X}.
\end{equation}
But, since $H = X^\mu \leq X^{1 - (k/2+1)\epsilon - 2\epsilon_1}$,
\begin{equation}
\label{4.6}
\begin{split}
H (X^{\mu + \epsilon})^{k/2+1} \log^{2k} {X} \leq & X^{1 -
(k/2+1)\epsilon - 2\epsilon_1} H^{k/2+1} X^{(k/2+1)\epsilon}
\log^{2k} X \\
\ll & H^{k/2+1} X^{1-\epsilon_1}.
\end{split}
\end{equation}
Combining (\ref{4.4}), (\ref{4.5}) and (\ref{4.6}), we have
\begin{equation}
\label{4.7}
\begin{split}
& \int_{0}^{H} \int_{1}^{X} \bigl(\psi(x+h) - \psi(x) -h \bigr)^k dx \;
dh \\
=& \nu_k H^{{k/2}+2} \int_{H/X}^{1/E} {1 \over \delta^2}
\Bigl(\log{1 \over E \delta}\Bigr)^{k/2} d\delta + \nu_k
X^{{k/2}+2} \int_{0}^{H/X} \delta^{k/2}
\Bigl(\log{1 \over E\delta}\Bigr)^{k/2} d\delta \\
&+ O_k(H^{{k/2}+1} X^{1-\epsilon_1})
\end{split}
\end{equation}
for $X^{2\epsilon} \leq H \leq X^{1-(k/2+1)\epsilon -
2\epsilon_1}$.

\bigskip

We now deduce (\ref{4.2}) from (\ref{4.7}). Set
$\eta=X^{-2\epsilon_1/3}$. For $X^{2\epsilon + \epsilon_1} \leq H
\leq X^{1-(k/2+1)\epsilon - 2\epsilon_1}/2$,
\begin{eqnarray*}
& &\int_{H}^{(1+\eta)H} \int_{1}^{X} (\psi(x+h) - \psi(x) - h)^k
dx \; dh \\
&=& \nu_k ((1+\eta)H)^{k/2+2} \int_{(1+\eta)H/X}^{1/E} {1 \over
\delta^2}
\Bigl(\log{1 \over E\delta}\Bigr)^{k/2} d\delta \\
& & - \nu_k H^{k/2+2} \int_{H/X}^{1/E} {1 \over \delta^2}
\Bigl(\log{1 \over E \delta}\Bigr)^{k/2} d\delta \\
& &+ \nu_k X^{k/2+2} \int_{H/X}^{(1+\eta)H/X} \delta^{k/2}
\Bigl(\log{1 \over E \delta}\Bigr)^{k/2} d\delta
+ O_k(H^{{k/2}+1} X^{1-\epsilon_1}) \\
&=& -\nu_k H^{k/2+2} \int_{H/X}^{(1+\eta)H/X} {1 \over \delta^2}
\Bigl(\log{1 \over E\delta}\Bigr)^{k/2} d\delta + \mu_k \eta
H^{k/2+2} \int_{H/X}^{1/E} {1 \over \delta^2} \Bigl(\log{1 \over
E\delta}\Bigr)^{k/2} d\delta \\
& &+ \nu_k X^{k/2+2} \int_{H/X}^{(1+\eta)H/X} \delta^{k/2}
\Bigl(\log{1 \over E \delta} \Bigr)^{k/2} d\delta \\
& &+ O_k\Bigl(\eta^2 H^{k/2+1} X (\log{X \over H})^{k/2}\Bigr) +
O_k(H^{k/2+1} X^{1-\epsilon_1}) \\
&=& \mu_k \eta H^{k/2+2} \int_{H/X}^{1/E} {1 \over \delta^2}
\Bigl(\log{1 \over E\delta}\Bigr)^{k/2} d\delta \\
& &+ O_k\Bigl(\eta^2 H^{k/2+1} X (\log{X \over H})^{k/2}\Bigr) +
O_k(H^{k/2+1} X^{1-\epsilon_1})
\end{eqnarray*}
by Lemma \ref{lemma2.1}. Therefore
\begin{equation}
\label{4.8}
\begin{split}
& \int_{H}^{(1+\eta)H} \int_{1}^{X} \bigl(\psi(x+h) - \psi(x) -h
\bigr)^k dx \; dh \\
=& \mu_k \eta H^{k/2+2} \int_{E}^{X/H} \bigl(\log{u \over E}
\bigr)^{k/2} du \\
&+ O_k\Bigl(\eta^2 H^{k/2+1} X (\log{X \over H})^{k/2}\Bigr) +
O_k(H^{{k/2}+1} X^{1-\epsilon_1}).
\end{split}
\end{equation}
Let $g(x,h)=f(x,H)$ for $H \leq h \leq (1+\eta)H$. Again, one can
check that $f(x,h)-g(x,h)=f(x+H,h-H)$. So,
\begin{equation}
\label{4.9}
\begin{split}
\int_{H}^{(1+\eta)H} \int_{1}^{X} \bigl(f(x,h)-g(x,h)\bigr)^k dx
dh &= \int_{0}^{\eta H} \int_{1+H}^{X+H} f(x,h)^k dxdh \\
&\ll_k \eta^{k/2+1} X H^{k/2+1} \Bigl(\log{X \over \eta
H}\Bigr)^{k/2}
\end{split}
\end{equation}
by (\ref{4.7}) as well as the choice of $\eta$ and the range of
$H$. Thus, by Lemma \ref{lemma2.3}, (\ref{4.8}) and (\ref{4.9}),
\begin{equation}
\label{4.10}
\begin{split}
\int_{H}^{(1+\eta)H} \int_{1}^{X} g(x,h)^k dx \; dh &\ll_k \int
\int |f(x,h)|^k + \int \int |f(x,h)-g(x,h)|^k \\
&\ll_k \eta X H^{k/2+1} (\log{X \over \eta H})^{k/2}+
H^{k/2+1}X^{1-\epsilon_1}.
\end{split}
\end{equation}
By Lemma \ref{lemma2.2} and Holder's inequality,
\begin{equation}
\label{4i1}
\begin{split}
&\int_{H}^{(1+\eta)H} \int_{1}^{X} f(x,h)^k - g(x,h)^k dx d\delta \\
=&\int_{H}^{(1+\eta)H} \int_{1}^{X} P(f,g)(f-g) +
\int_{H}^{(1+\eta)H} \int_{1}
^{X} (f-g)^k \\
\ll& \Bigl(\int \int |P(f,g)|^{k/(k-1)}\Bigr)^{(k-1)/k} \Bigl(\int
\int |f-g|^k\Bigr)^{1/k} + \int \int |f-g|^k \\
=&K_1^{(k-1)/k} K_2^{1/k} + K_2
\end{split}
\end{equation}
where $P(x,y)$ is a homogeneous polynomial of degree $k-1$. From
(\ref{4.9}),
\begin{equation}
\label{4i2} K_2 \ll_k \eta^{k/2+1} X H^{k/2+1} \Bigl(\log{X \over
\eta H}\Bigr)^{k/2}.
\end{equation}
And similar to the proof in Theorem \ref{theorem3.1},
\begin{equation}
\label{4i3} K_1 \ll_k \int \int f^k+g^k \ll_k \eta X
H^{k/2+1}(\log{X \over \eta H})^{k/2}
\end{equation}
by (\ref{4.8}) and (\ref{4.10}). Consequently, by (\ref{4i1}),
(\ref{4i2}), (\ref{4i3}) and Lemma \ref{lemma2.3},
\begin{equation}
\label{4.11} \int \int f^k -g^k \ll_k \eta^{3/2} X H^{k/2+1}
\Bigl(\log{X \over \eta H} \Bigr)^{k/2}.
\end{equation}
Therefore, by (\ref{4.11}) and (\ref{4.8}),
\begin{eqnarray*}
& &\eta H \int_{1}^{X} \bigl(\psi(x+h) - \psi(x) -h \bigr)^k dx =
\int \int g^k \\
&=&\int \int f^k + O_k\Bigl(X H^{k/2+1} \eta^{3/2} (\log{X \over
\eta H})^{k/2}\Bigr) \\
&=& \mu_k \eta H^{k/2+2} \int_{E}^{X/H} \bigl(\log{u \over
E}\bigr)^{k/2} du + O_k(H^{{k/2}+1} X^{1-\epsilon_1}) \\
& &+ O_k\Bigl(X H^{k/2+1} \eta^{3/2} (\log{X \over \eta
H})^{k/2}\Bigr).
\end{eqnarray*}
Divide through by $\eta H$ and recall $\eta = X^{-2\epsilon_1/3}$,
we get the theorem with $\epsilon_2 = {\epsilon_1 \over 4}$.
\section{Numerical evidence}

In Montgomery and Soundararajan [\ref{MS1}], they got some
numerical data for the actual values of $M_k(X;h)$. One has the
following table:

\bigskip

For $X = 10^{10}$ and $h = 10^5$.

\bigskip

\begin{math}
\begin{array}{lll}
k & \mbox{Actual value of }M_k(X;h) &\mbox{Result from
formula (i) of Theorem \ref{theorem1.1}} \\
2 & 9.0663\ast 10^{15} & 9.0978\ast 10^{15} \\
4 & 2.4995\ast 10^{22} & 2.5131\ast 10^{22} \\
6 & 1.1573\ast 10^{29} & 1.1675\ast 10^{29}
\end{array}
\end{math}

\bigskip

Using a C program, we get some numerical evidence in support of
the truth of (ii) in Theorem \ref{theorem1.1}.

\bigskip

For $X=10^8$ and $\delta=10^{-4}$:

\bigskip

\begin{math}
\begin{array}{lll}
k & \mbox{Actual value of }\widetilde{M}_k(X;\delta) &
\mbox{Result from formula (ii) of Theorem \ref{theorem1.1}} \\
2 & 4.0075\ast 10^{12} & 3.8976\ast 10^{12} \\
4 & 6.5161\ast 10^{17} & 6.0766\ast 10^{17} \\
6 & 1.9592\ast 10^{23} & 1.7763\ast 10^{23}
\end{array}
\end{math}

\bigskip
For $X=10^{10}$ and $\delta=10^{-5}$:
\bigskip

\begin{math}
\begin{array}{lll}
k & \mbox{Actual value of }\widetilde{M}_k(X;\delta) &
\mbox{Result from formula (ii) of Theorem \ref{theorem1.1}} \\
2 & 5.0527\ast 10^{15} & 5.0485\ast 10^{15} \\
4 & 1.0210\ast 10^{22} & 1.0195\ast 10^{22} \\
6 & 3.8645\ast 10^{28} & 3.8602\ast 10^{28}
\end{array}
\end{math}


Tsz Ho Chan\\
American Institute of Mathematics\\
360 Portage Avenue\\
Palo Alto, CA 94306\\
USA\\
thchan@aimath.org

\end{document}